\documentclass{proc-l}

\newtheorem{theorem}{Theorem}
\newtheorem{lemma}{Lemma}

\theoremstyle{remark}
\newtheorem{remark}{Remark}

\newcommand{\free}{{\mathcal F}}

\newcommand{\trans}{\top}

\begin{document}

\title
{Non-commutative positive kernels and their matrix evaluations}

\author[D. S. Kalyuzhny\u{\i}-Verbovetzki\u{\i}]{Dmitry S. Kalyuzhny\u{\i}-Verbovetzki\u{\i}}
\address{Department of Mathematics, Ben-Gurion University of the Negev, Beer Sheva, Israel 84105}
\email{dmitryk@math.bgu.ac.il}
\thanks{The first author was supported by the Center for Advanced Studies in Mathematics, Ben-Gurion University of the Negev.}

\author[V. Vinnikov]{Victor Vinnikov}
\address{Department of Mathematics, Ben-Gurion University of the Negev, Beer Sheva, Israel 84105}
\email{vinnikov@math.bgu.ac.il}
\thanks{The second author was partially supported by the Israel Science Foundation
Grant 322/00-1.}





\begin{abstract}
We show that a formal power series in $2N$ non-commuting indeterminates
is a positive non-commutative kernel if and only if
the kernel on $N$-tuples of matrices of any size
obtained from this series by matrix substitution is positive. 
We present  two versions of this result related to  different classes of matrix substitutions. 
In the general case we consider substitutions of jointly nilpotent $N$-tuples of matrices, 
and thus the question of convergence does not arise. 
In the ``convergent'' case we consider substitutions of $N$-tuples of matrices 
from a neighborhood of zero where the series converges. Moreover, in the first case 
the result can be improved: the positivity of a non-commutative kernel is guaranteed by 
the positivity of its values on the diagonal, i.e., on pairs of coinciding 
jointly nilpotent $N$-tuples of matrices. 
In particular this yields an analogue
of a recent result of Helton on non-commutative sums-of-squares representations
for the class of hereditary non-commutative polynomials.
We show by an example that the improved formulation
does not apply in the ``convergent'' case.
\end{abstract}

\maketitle

\section{Introduction and statement of the results}
\label{s:intro}

Several recent results show a relationship between properties 
of a polynomial or of a formal power series
in say $N$ non-commuting indeterminates
and the properties of the matrix-valued functions on $N$-tuples of matrices of all sizes 
obtained from this polynomial or this formal power series by matrix substitution; e.g.,
see \cite{Hel,McC,HelMcCPut,BMG,AKV}.

The objective of this paper is to establish a result of this type 
in the framework of positive kernels \cite{Aro,Sai,ADRdS}.
We recall that if $\Omega$ is any set and ${\mathcal E}$ is a Hilbert space,
an operator valued function
$K \colon \Omega \times \Omega \to {\mathcal L}({\mathcal E})$\footnote{
${\mathcal L}({\mathcal E})$ denotes the $C^*$-algebra of bounded linear operators
on ${\mathcal E}$.}
is called a positive ${\mathcal L}({\mathcal E})$-valued kernel on $\Omega$
if for any integer $l$,
any points $\omega^{(1)},\ldots,\omega^{(l)} \in \Omega$ and 
any vectors $h^{(1)},\ldots,h^{(l)}\in\mathcal{E}$, one has
$$
\sum_{j,k=1}^l\left\langle 
K(\omega^{(j)},\omega^{(k)})h^{(k)},h^{(j)}\right\rangle_{\mathcal{E}}\geq 0.
$$
Alternatively, there exists an auxilliary Hilbert space ${\mathcal H}$ and an operator
valued function $H \colon \Omega \to {\mathcal L}({\mathcal H},{\mathcal E})$\footnote{
${\mathcal L}({\mathcal H},{\mathcal E})$ denotes the Banach space of bounded linear
operators from ${\mathcal H}$ to ${\mathcal E}$.}
so that the factorization
$$
K(\omega,\omega') = H(\omega) H(\omega')^*
$$
holds for all $\omega,\omega' \in \Omega$.

Let $\free_N$ denote the free semigroup with $N$ generators  $g_1,\ldots,g_N$ and unit element 
$\emptyset$. For  $w=g_{j_1}\cdots g_{j_m}$ and $w'=g_{k_1}\cdots g_{k_l}$ from 
$\free_N$ (words in the alphabet  $\{ g_1,\ldots,g_N\} )$ 
the semigroup multiplication is the concatenation:
$$ww'=g_{j_1}\cdots g_{j_m}g_{k_1}\cdots g_{k_l},$$
and $\emptyset\in\free_N$ plays a role of the empty word: for any $w\in\free_N$,
$$\emptyset w=w\emptyset =w.$$
For a word  $w=g_{j_1}\cdots g_{j_m}\in\free_N$ define its transpose 
$w^\trans=g_{j_m}\cdots g_{j_1}\in\free_N$ and its length 
$|w|=m\in\mathbb{N}$, and define $\emptyset^\trans=\emptyset,\ |\emptyset |=0$.

Let $z_1,\ldots,z_d$ be non-commuting indeterminates, and $w=g_{j_1}\cdots g_{j_m}\in\free_d$. Set
$$z^w=z_{j_1}\cdots z_{j_m},\quad z^\emptyset =1.$$
For any vector space ${\mathcal L}$ we denote by
${\mathcal L}\langle\langle
z_1,\ldots,z_d \rangle\rangle$ the space of all formal power series
in $z_1,\ldots,z_d$ with coefficients in ${\mathcal L}$.

A formal power series in $2N$ non-commuting indeterminates $z=(z_1,\ldots,z_N)$
and $z'=(z'_1,\ldots,z'_N)$,
\begin{equation} \label{fps}
K(z,z')=\sum_{w,w'\in{\mathcal F}_N}K_{w,w'}z^wz^{\prime w^{\prime \trans}} \in
{\mathcal L}({\mathcal E})\langle\langle
z, z' \rangle\rangle,
\end{equation}
is called a {\em positive non-commutative kernel} if
$(w,w') \mapsto K_{w,w'}$ is a positive ${\mathcal L}({\mathcal E})$-valued
kernel on ${\mathcal F}_{N}$, i.e.,
for any integer $l$, any words $w^{(1)},\ldots,w^{(l)}\in\free_N$ and any
vectors $h^{(1)},\ldots,h^{(l)}\in\mathcal{E}$, one has
$$\sum_{j,k=1}^l\left\langle K_{w^{(j)},w^{(k)}}h^{(k)},h^{(j)}\right\rangle_{\mathcal{E}}\geq 0.$$
Alternatively,
there is an auxiliary Hilbert space ${\mathcal H}$ and a
non-commutative formal power series 
\begin{equation} \label{zps}
H(z) = \sum_{w \in \free_N} H_w z^w \in {\mathcal
L}({\mathcal H}, {\mathcal E})\langle\langle z \rangle\rangle,
\end{equation}
such that
$$
K(z,  z') = H(z) H(z')^{*}
$$
where we use the convention
$$ (z^{\prime w})^{*} = z^{\prime w^{\trans}}.  $$
This notion of positivity was studied extensively in \cite{FRKHS}.\footnote{
In \cite{FRKHS}, it is assumed that $z_i$ commutes with $z_j'$.
For our purposes it is more convenient to deal with $2N$ non-commuting indeterminates;
notice that it does not really matter since in all our formulae every $z_j'$
appears to the right of every $z_i$, in every monomial.}

Let $Z=(Z_1,\ldots,Z_N)\in (\mathbb{C}^{n\times n})^N$, i.e., $Z$ is an $N$-tuple of $n\times n$ matrices with complex entries, and  $w=g_{j_1}\cdots g_{j_m}\in\free_N$. Set
$$Z^w=Z_{j_1}\cdots Z_{j_m},\quad Z^\emptyset =I_n.$$
Let $n,r\in\mathbb{N}$. Define the set $\mbox{Nilp}_N(n,r)$ of $N$-tuples $Z$ of $n\times n$ matrices with complex entries, which are \emph{jointly nilpotent of rank at most} $r$, i.e., $Z^w=0$ for all $w\in\free_N$ such that $|w|\geq r$. 
Define also the set $\mbox{Nilp}_N(n)=\bigcup_{r=1}^\infty\mbox{Nilp}_N(n,r)$ of 
\emph{jointly nilpotent} $N$-tuples of $n\times n$ matrices\footnote{
Let us remark that a jointly nilpotent tuple of matrices is jointly similar 
to a tuple of strictly upper triangular matrices 
--- the authors thank Leonid Gurvits 
for pointing this out.}. 

\begin{theorem}\label{main}
A formal power series $K(z,z') \in {\mathcal L}({\mathcal E})\langle\langle z, z' \rangle\rangle$ as in \eqref{fps}
is a positive non-commutative kernel if and only if \begin{equation} \label{mps}
K(Z,Z') := \sum_{w,w'\in\free_N}K_{w,w'}\otimes Z^w{{Z'}^*}^{{w'}^\trans} \in
{\mathcal L}(\mathcal{E}\otimes\mathbb{C}^n)
\end{equation}
where $Z,Z' \in \rm{Nilp}_N(n)$, 
${Z'}^*=({Z'}^*_1,\ldots,{Z'}^*_N)$
(notice that the sum is therefore finite!),  
is a positive ${\mathcal L}(\mathcal{E}\otimes\mathbb{C}^n)$-valued
kernel on $\rm{Nilp}_N(n)$ for every $n\in\mathbb{N}$, i.e., for any integer $l$, any jointly nilpotent $N$-tuples of matrices 
$Z^{(1)},\ldots,Z^{(l)}\in\rm{Nilp}_N(n)$ and any
vectors $x^{(1)},\ldots,x^{(l)}\in\mathcal{E}\otimes\mathbb{C}^n$, one has
$$\sum_{j,k=1}^l\left\langle K(Z^{(j)},Z^{(k)})x^{(k)},x^{(j)}\right\rangle
_{\mathcal{E}\otimes\mathbb{C}^n}\geq 0.$$ 
\end{theorem}

One direction of the theorem is obvious.
If $K(z,z')$ is a positive non-commutative kernel then
we have a factorization
$K(z,z') = H(z) H(z')^*$
for some formal power series 
$H(z)$ as in \eqref{zps}
and some Hilbert space ${\mathcal H}$. Therefore for any
$Z, Z' \in \mbox{Nilp}_N(n)$ we have
$K(Z,Z') = H(Z) H(Z')^*,$
where
\begin{equation} \label{zpssub}
H(Z) := \sum_{w \in \free_N} H_w \otimes Z^w \in {\mathcal L}
({\mathcal H} \otimes {\mathbb C}^n,{\mathcal E} \otimes {\mathbb C}^n)
\end{equation}
(notice that the sum is finite).
Therefore for each $n\in\mathbb{N}$ we have that $K(Z,Z')$ is a positive ${\mathcal L}(\mathcal{E}\otimes\mathbb{C}^n)$-valued
kernel on $\mbox{Nilp}_N(n)$.

Section \ref{s:proof1} contains the proof of the other direction of the theorem.

Notice that by using nilpotent matrix
substitutions one avoids any convergence assumptions on the formal power series.
There is also a ``convergent'' version of Theorem \ref{main}.
A formal power series
$K(z,z') \in {\mathcal L}({\mathcal E})\langle \langle z, z' \rangle\rangle$
as in \eqref{fps} is called {\em convergent} if
for every $n\in\mathbb{N}$ there exists a connected open neighbourhood $U_n$ of $0$ in
$(\mathbb{C}^{n\times n})^N$ such that the series $K(Z,Z')$ as in \eqref{mps}
converges uniformly on compact subsets of $U_n\times U_n$ in the norm of ${\mathcal L}({\mathcal E} \otimes {\mathbb C}^n)$ (we assume that $\mathcal{F}_N$ and $\mathcal{F}_{2N}\cong\mathcal{F}_N\times\mathcal{F}_N$ are ordered, say, lexicographically).

\begin{theorem}\label{conv}
A convergent formal power series $K(z,z') \in {\mathcal L}({\mathcal E})\langle\langle z, z' \rangle\rangle$ (with respect to a set of neighbourhoods $U_n$ of $0$ in $(\mathbb{C}^{n\times n})^N,\ n\in\mathbb{N}$)
is a positive non-commutative kernel if and only if  $K(Z,Z')$
 is a positive ${\mathcal L}(\mathcal{E}\otimes\mathbb{C}^n)$-valued kernel on $U_n$ for every $n\in\mathbb{N}$.
\end{theorem}

We prove this theorem in Section~\ref{s:proof2}.

In the ``nilpotent'' case there is a considerable strengthening of Theorem \ref{main}.
The authors thank the anonymous referee for suggesting this result and its proof.
\begin{theorem}\label{impr}
A formal power series $K(z,z') \in {\mathcal L}({\mathcal E})\langle\langle z, z' \rangle\rangle$ as in
\eqref{fps}
is a positive non-commutative kernel if and only if for every $n\in\mathbb{N}$ and $Z\in {\rm Nilp}_N(n)$,
\begin{equation} \label{diag}
K(Z,Z) := \sum_{w,w'\in\free_N}K_{w,w'}\otimes Z^w{Z^*}^{{w'}^\trans} \in
{\mathcal L}(\mathcal{E}\otimes\mathbb{C}^n)
\end{equation}
is a positive semidefinite operator, i.e., for any $x\in\mathcal{E}\otimes\mathbb{C}^n$ one has
$$\left\langle K(Z,Z)x,x\right\rangle_{\mathcal{E}\otimes\mathbb{C}^n}\geq 0.$$ 
\end{theorem}

One direction of the theorem is obvious due to the remark following the formulation of 
Theorem~\ref{main} and the fact that the operator $H(Z)H(Z)^*$ is positive semidefinite. 
Section~\ref{s:proofs3,4} contains the proof of the other direction of Theorem \ref{impr}.
While Theorem~\ref{impr} implies Theorem~\ref{main}, Section~\ref{s:proof1} contains
an independent proof of Theorem~\ref{main} since the argument of this
proof is also used in  the proof of Theorem~\ref{conv}. 
Let us remark also that the statement analogous to the one of Theorem~\ref{impr} 
does not hold in the ``convergent'' case. 
Indeed, set $K(z,z'):=1-zz'\in\mathbb{C}\left\langle z\right\rangle$ (here $N=1$). 
Clearly, $K(Z,Z)>0$ for a matrix $Z$ (of any size) close to the origin, 
however $K(z,z')$ is not a positive kernel since its matrix of coefficients 
{\scriptsize $\left(\begin{array}{cc}
 1 & 0\\
 0 & -1
\end{array}\right)$} is indefinite.
 
Theorem \ref{impr} implies a factorization result for a class of non-commutative polynomials.
For a vector space $\mathcal{L}$ we denote by 
$\mathcal{L}\left\langle z_1,\ldots,z_d\right\rangle$ 
the space of all polynomials in non-commuting indeterminates 
$z_1,\ldots,z_d$ with coefficients in $\mathcal{L}$, i.e., 
a subspace in $\mathcal{L}\left\langle\left\langle z_1,\ldots,z_d\right\rangle\right\rangle$ 
consisting of formal power series with a finitely supported set of coefficients. 
An operator-valued polynomial in $2N$ non-commuting indeterminates $z=(z_1,\ldots, z_N)$ 
and $z'=(z_1',\ldots,z_N')$ is called called {\em hereditary} if it is of the form
\begin{equation} \label{heredpoly}
K(z,z')=\sum_{w,w'\in\mathcal{F}_N:\,|w|\leq m,|w'|\leq m}
K_{w,w'}z^w{z'}^{{w'}^T}\in\mathcal{L(E)}\left\langle z,z'\right\rangle,
\end{equation}
where $\mathcal{E}$ is a Hilbert space, i.e., every $z_j'$ appears to the right of every
$z_i$, in every monomial.
\begin{theorem}\label{fact}
A hereditary polynomial $K(z,z') \in \mathcal{L(E)}\left\langle z,z'\right\rangle$
as in \eqref{heredpoly} satisfies $K(Z,Z) \geq 0$ for all $Z \in {\rm Nilp}_N(n)$ ($n=1,2,\ldots$)
if and only if there exist an auxiliary Hilbert space $\mathcal{H}$ 
(of dimension at most $\dim(\mathcal{E})\sum_{j=0}^mN^j$) and a polynomial 
$$H(z)=\sum_{w\in\mathcal{F}_N:\,|w|\leq m}H_wz^w\in\mathcal{L}(\mathcal{H,E})\left\langle z\right\rangle$$ 
such that 
$$K(z,z')=H(z)H(z')^*.$$
Furthermore, for this factorization to hold it is enough to assume
that $K(Z,Z) \geq 0$ holds for all
$Z \in {\rm Nilp}_N(\sum_{j=0}^mN^j,m+1)$.
\end{theorem}
We deduce this theorem in Section~\ref{s:proofs3,4}. 
Theorem~\ref{fact} is an analogue of the main result of \cite{Hel} 
(for related results and generalizations see \cite{McC} and \cite{HelMcCPut}),
for hereditary non-commutative polynomials and 
jointly nilpotent test matrices. 
 
In view of the results of this paper it is of interest to investigate  the relationship between the non-commutative formal
reproducing kernel Hilbert space associated to a positive non-commutative kernel (such spaces are studied in \cite{FRKHS})
and the reproducing kernel Hilbert spaces associated
to the kernels obtained from the non-commutative kernel by matrix substitutions.

\section{Proof of Theorem~\ref{main}}
\label{s:proof1}

For the proof of Theorem~\ref{main} we need the following lemma.

\begin{lemma}\label{lem}
Let $\mathcal{U}$ be a Hilbert space and let
an ${\mathcal L}(\mathcal{U})$-valued polynomial function 
$$
P(\lambda,\lambda')=
\sum\limits_{t,t'\in\mathbb{Z}^N_+:\ 0\leq |t|,|t'|\leq m}
P_{t,t'}\lambda^t(\overline{\lambda'})^{t'}
$$
be a positive kernel on ${\mathbb C}^N$. 
Here for $\lambda=(\lambda_1,\ldots,\lambda_N)\in\mathbb{C}^N$ and 
$t=(t_1,\ldots,t_N)\in\mathbb{Z}^N_+$ 
we set $\lambda^t={\lambda_1}^{t_1}\cdots\lambda_N^{t_N}$ and 
$|t|=t_1+\cdots +t_N$. 
Then the matrix $M_P=(P_{t,t'})_{0\leq |t|,|t'|\leq m}$ 
(for a certain order on $\mathbb{Z}^N_+$) 
defines a positive semidefinite operator on $\mathcal{U}^{\frac{(m+N)!}{m!N!}}$.
\end{lemma}

This holds also for sesquianalytic rather than polynomial kernels
(of course in this case one has to consider the matrices
$(P_{t,t'})_{0\leq |t|,|t'|\leq m}$ for any integer $m$)
and the converse is true as well; notice that these facts are
a ``commutative convergent'' analogue of Theorem~\ref{main}.
They are well known for $N=1$ (see, e.g., \cite[Lemma 1.1.5]{ADRdS})
and the proofs in the general case are analogous.
We provide a proof of the lemma for the sake of completeness.

\begin{proof}[Proof of Lemma \ref{lem}]
For an arbitrary $u=(u_t)_{0\leq |t|\leq m}\in\mathcal{U}^{\frac{(m+N)!}{m!N!}}$ we have
\begin{multline*}
\left\langle M_Pu,u\right\rangle_{\mathcal{U}^{\frac{(m+N)!}{m!N!}}} = 
\sum_{t,t'\in\mathbb{Z}^N_+:\ 0\leq |t|,|t'|\leq m}\left\langle  
P_{t,t'}u_{t'},u_t\right\rangle_\mathcal{U}\\ =
\sum_{t,t'\in\mathbb{Z}^N_+:\ 0\leq |t|,|t'|\leq m}\  
\int\limits_{\mathbb{T}^N\times\mathbb{T}^N}\left\langle  
P(\lambda,\lambda')u_{t'},u_t\right\rangle_\mathcal{U}{\overline{\lambda}}^t{\lambda'}^{t'}\,dm_N(\lambda)dm_N(\lambda'),
\end{multline*}
where $dm_N(\cdot)$ is the normalized Lebesgue measure on $\mathbb{T}^N$. 
For an arbitrary $n\in\mathbb{N}$, 
let $\{\Delta_j\}_{j=1}^n$ be a partition of $\mathbb{T}^N$ such that 
$m_N(\Delta_j)=\frac{1}{n}$ and let $\lambda^{(j)}\in\Delta_j$, $j=1,\ldots,n$, 
be arbitrary points. Then
\begin{multline*}
\left\langle M_Pu,u\right\rangle_{\mathcal{U}^{\frac{(m+N)!}{m!N!}}} \\ = 
\sum_{0\leq |t|,|t'|\leq m}\lim_{n\to\infty}\frac{1}{n^2}
\sum_{j,k=1}^{n}\left\langle  
P(\lambda^{(j)},\lambda^{(k)})
u_{t'},u_t\right\rangle_\mathcal{U}{(\overline{\lambda^{(j)}})}^t{(\lambda^{(k)})}^{t'} \\=
\lim_{n\to\infty}\frac{1}{n^2}\sum_{j,k=1}^{n}\left\langle  
P(\lambda^{(j)},\lambda^{(k)})
\sum_{0\leq |t'|\leq m}{(\lambda^{(k)})}^{t'}u_{t'},
\sum_{0\leq |t|\leq m}{(\lambda^{(j)})}^tu_t\right\rangle_\mathcal{U}\geq 0,
\end{multline*}
as desired.
\end{proof}

\begin{proof}[Proof of Theorem~\ref{main}.] As we already said in Section~\ref{s:intro}, one direction of this theorem is obvious, and thus only the other one is left to prove. 
Assume that a formal power series 
$K(z,z') \in {\mathcal L}({\mathcal E}) \langle\langle z, z' \rangle\rangle$ as in \eqref{fps} is such that the function $K(Z,Z')$ as in \eqref{mps} is a positive ${\mathcal L}(\mathcal{E}\otimes\mathbb{C}^n)$-valued
kernel on $\rm{Nilp}_N(n)$ for every $n\in\mathbb{N}$.

Let  $m\in\mathbb{N}$, $m>1$. Set $n=(N+1)^m$. 
Define a linear operator $S$ on the space 
$\mathbb{C}^n\cong(\mathbb{C}^{N+1})^{\otimes m}$ by its action on the basis vectors 
$e_{i_1}\otimes\cdots\otimes e_{i_m}$ where $e_j$, $j=1,\ldots,N+1$ 
are the standard basis vectors in $\mathbb{C}^{N+1}$:
$$
S(e_{i_1}\otimes\cdots\otimes e_{i_m})=  e_{i_m}\otimes e_{i_1}\otimes\cdots\otimes e_{i_{m-1}}.
$$
The operator $S$ is represented 
in the basis $e_{i_1}\otimes\cdots\otimes e_{i_m},\ i_1,\ldots,i_m\in\{ 1,\ldots,N+1\}$,
by a $n\times n$ permutation matrix.
Let $E_{i,j}\in\mathbb{C}^{(N+1)\times (N+1)}$, $i,j=1,\ldots,N+1$, be matrices of the form
$$
(E_{i,j})_{\mu,\nu}=\left\{
\begin{array}{cc}
1,&(\mu,\nu)=(i,j),\\0, & (\mu,\nu)\neq (i,j).
\end{array}
\right.
$$
Define 
$$
Z_k=\left(E_{k+1,1}\otimes I_{N+1}^{\otimes m-1}\right)S\in\mathbb{C}^{n\times n}
\cong\left(\mathbb{C}^{(N+1)\times (N+1)}\right)^{\otimes m},\ k=1,\ldots,N.
$$
For $\lambda=(\lambda_1,\ldots,\lambda_N)\in\mathbb{C}^N$, set
$$
Z(\lambda)=(\lambda_1Z_1,\ldots,\lambda_NZ_N)\in\left(\mathbb{C}^{n\times n}\right)^N.
$$
For $w=g_{j_1}\cdots g_{j_{|w|}}\in\free_N:\ 0<|w|\leq m$, one has
$$
Z(\lambda)^w=\lambda_{j_1}Z_{j_1}\cdots\lambda_{j_{|w|}}Z_{j_{|w|}}=
\lambda^w(E_{j_1+1,1}\otimes\cdots\otimes E_{j_{|w|}+1,1}\otimes I_{N+1}^{\otimes m-|w|})S^{|w|}
$$
(in the case $|w|=m$, 
the last term in the tensor product above disappears, and $S^{|w|}=I_n$), 
or, on the basis vectors,
$$
Z(\lambda)^w(e_{i_1}\otimes\cdots\otimes e_{i_m})=
\lambda^w(E_{j_1+1,1}e_{i_{m-|w|+1}}\otimes\cdots\otimes E_{j_{|w|}+1,1}e_{i_m}
\otimes e_{i_1}\otimes\cdots\otimes e_{i_{m-|w|}}).
$$

Let us remark that since $\lambda_k$, $k=1,\ldots,N$, are scalars, and thus commute, 
we may rewrite $\lambda^w$ as $\lambda^{t(w)}$, 
where $t:\free_N\to\mathbb{Z}^N_+$ is the \emph{abelianization map}:
$t(w)=(t_1(w),\ldots,t_N(w))$ with the non-negative integer number $t_k(w)$ 
equal to the numbers of times that the letter $g_k$ appears in the word $w$. 
Thus, we have 
$$
\lambda^w=\lambda^{t(w)}=\lambda_1^{t_1(w)}\cdots\lambda_N^{t_N(w)}.
$$

For $w=g_{j_1}\cdots g_{j_{|w|}}\in\free_N:\ |w|= m+1$, one has 
\begin{align*}
Z(\lambda)^w 
&= \lambda^{t(w)}(E_{j_1+1,1}\otimes I_{N+1}^{\otimes m-1})S(E_{j_2+1,1}\otimes\cdots\otimes E_{j_{m+1},1})  \\ 
&= \lambda^{t(w)}(E_{j_1+1,1}E_{j_{m+1},1}\otimes E_{j_2+1,1}\otimes\cdots\otimes E_{j_m,1})S=0,
\end{align*}
since $E_{j+1,1}E_{k+1,1}=0$ for any $j,k\in\{ 1,\ldots,N\}$. Thus, $Z(\lambda)\in\mbox{Nilp}_N(n,m+1)$, and
by the assumption, 
\begin{align*}
P(\lambda,\lambda') 
&= K(Z(\lambda),Z(\lambda'))=\sum_{w,w'\in\free_N:\ 0\leq |w|,|w'|\leq m}K_{w,w'}\otimes Z(\lambda)^w{Z(\lambda')^*}^{{w'}^\trans}\\
&= \sum_{t,t'\in\mathbb{Z}^N_+:\ 0\leq |t|,|t'|\leq m}\left(\sum_{w,w'\in\free_N:\ t(w)=t,t(w')=t'}K_{w,w'}\otimes Z^w{Z^*}^{{w'}^\trans}\right)\lambda^t{(\overline{\lambda'})}^{t'}
\end{align*}
is a polynomial ${\mathcal L}(\mathcal{E}\otimes\mathbb{C}^N)$-valued function which is a positive 
kernel on ${\mathbb C}^N$. 
By Lemma~\ref{lem}, the coefficients of this kernel serve as the operator blocks $P_{t,t'}\in 
{\mathcal L}(\mathcal{E}\otimes\mathbb{C}^n)$ of the matrix $M_P=(P_{t,t'})_{0\leq |t|,|t'|\leq m}$ (for some order on $\mathbb{Z}^N_+$) which represents a positive semidefinite operator on the Hilbert space $(\mathcal{E}\otimes\mathbb{C}^n)^{\frac{(m+N)!}{m!N!}}$. Therefore, for  arbitrary vectors
$$u_t=\sum_{1\leq i_1,\ldots,i_m\leq N+1}u^{(t)}_{i_1,\ldots,i_m}\otimes e_{i_1}\otimes\cdots\otimes e_{i_m}\in\mathcal{E}\otimes\mathbb{C}^n\cong\mathcal{E}\otimes(\mathbb{C}^{N+1})^{\otimes m},$$
where $t\in\mathbb{Z}^N_+:\ 0\leq |t|\leq m,$ we can write down:
 \begin{align*}
0 & \leq  \sum_{t,t'\in\mathbb{Z}^N_+:\ 0\leq |t|,|t'|\leq m}\left\langle P_{t,t'}u_{t'},u_t\right\rangle_{\mathcal{E}\otimes\mathbb{C}^n}\\
&= \left(\sum_{|t|=|t'|}+\sum_{|t|>|t'|}+\sum_{|t|<|t'|}\right)\left\langle P_{t,t'}u_{t'},u_t\right\rangle_{\mathcal{E}\otimes\mathbb{C}^n} = A+B+C.
 \end{align*}
Let us calculate each of the sums $A$, $B$ and $C$ separately:
\allowdisplaybreaks
 \begin{align*}
A &= \sum_{|t|=|t'|}\left\langle P_{t,t'}u_{t'},u_t\right\rangle_{\mathcal{E}\otimes\mathbb{C}^n} \\ 
&= \sum_{|t|=|t'|}\left\langle \left(\sum_{w,w'\in\free_N:\ t(w)=t,\, t(w')=t'}K_{w,w'}\otimes Z^w{Z^*}^{{w'}^\trans}\right)u_{t'},u_t\right\rangle_{\mathcal{E}\otimes\mathbb{C}^n} \\ 
&= \sum_{|t|=|t'|}\, \sum_{t(w)=t,\, t(w')=t'}\left\langle  K_{w,w'}\otimes \left(E_{j_1+1,1}\otimes\cdots\otimes E_{j_{|w|}+1,1}\otimes I_{N+1}^{\otimes m-|w|}\right)S^{|w|} \right.\\
 &\times   \left. S^{*|w'|}\left(E_{1,k_1+1}\otimes\cdots\otimes E_{1,k_{|w'|}+1}\otimes I_{N+1}^{\otimes m-|w'|}\right) u_{t'},u_t\right\rangle_{\mathcal{E}\otimes\mathbb{C}^n} \\
 &= \sum_{|t|=|t'|}\, \sum_{t(w)=t,\, t(w')=t'}\left\langle  K_{w,w'}\otimes E_{j_1+1,k_1+1}\otimes\cdots\otimes E_{j_{|w|}+1,k_{|w'|}+1}\otimes I_{N+1}^{\otimes m-|t|}\right. \\
 &\times   \sum_{1\leq i_1,\ldots,i_m\leq N+1}u^{(t')}_{i_1,\ldots,i_m}\otimes e_{i_1}\otimes\cdots\otimes e_{i_m},\\
 &  \left. \sum_{1\leq s_1,\ldots,s_m\leq N+1}u^{(t)}_{s_1,\ldots,s_m}\otimes e_{s_1}\otimes\cdots\otimes e_{s_m}\right\rangle_{\mathcal{E}\otimes\mathbb{C}^n} \\
 &=  \sum_{|t|=|t'|}\, \sum_{t(w)=t,\, t(w')=t'}\, \sum_{1\leq i_{|t'|+1},\ldots,i_m\leq N+1}\left\langle  K_{w,w'}\right.
  u^{(t')}_{k_1+1,\ldots,k_{|w'|}+1,i_{|t'|+1},\ldots,i_m}, \\ 
  &  \left. u^{(t)}_{j_1+1,\ldots,j_{|w|}+1,i_{|t'|+1},\ldots,i_m} \right\rangle_{\mathcal{E}} ;
 \end{align*}
 \begin{align*}
 B &= \sum_{|t|>|t'|}\left\langle P_{t,t'}u_{t'},u_t\right\rangle_{\mathcal{E}\otimes\mathbb{C}^n}  \\
 &= \sum_{|t|>|t'|}\left\langle \left(\sum_{w,w'\in\free_N:\ t(w)=t,\, t(w')=t'}K_{w,w'}\otimes Z^w{Z^*}^{{w'}^\trans}\right)u_{t'},u_t\right\rangle_{\mathcal{E}\otimes\mathbb{C}^n} \\ 
&= \sum_{|t|>|t'|}\, \sum_{t(w)=t,\, t(w')=t'}\left\langle  K_{w,w'}\otimes \left(E_{j_1+1,1}\otimes\cdots\otimes E_{j_{|w|}+1,1}\otimes I_{N+1}^{\otimes m-|w|}\right)S^{|w|} \right.\\
 &\times   \left. S^{*|w'|}\left(E_{1,k_1+1}\otimes\cdots\otimes E_{1,k_{|w'|}+1}\otimes I_{N+1}^{\otimes m-|w'|}\right) u_{t'},u_t\right\rangle_{\mathcal{E}\otimes\mathbb{C}^n} \\
 &= \sum_{|t|>|t'|}\, \sum_{t(w)=t,\, t(w')=t'}\left\langle  K_{w,w'}\otimes E_{j_1+1,1}\otimes\cdots\otimes E_{j_{|w|-|w'|}+1,1}\otimes\right. \\ 
 &\otimes E_{j_{|w|-|w'|+1}+1,k_1+1}\otimes\cdots\otimes E_{j_{|w|}+1,k_{|w'|}+1}\otimes I_{N+1}^{\otimes m-|t|} \\
 &\times \sum_{1\leq i_1,\ldots,i_m\leq N+1}u^{(t')}_{i_1,\ldots,i_m}\otimes e_{i_{m-|w|+|w'|+1}}\otimes\cdots\otimes e_{i_m}\otimes e_{i_1}\otimes \ldots\otimes e_{i_{m-|w|+|w'|}},\\
 &  \left. \sum_{1\leq s_1,\ldots,s_m\leq N+1}u^{(t)}_{s_1,\ldots,s_m}\otimes e_{s_1}\otimes\cdots\otimes e_{s_m}\right\rangle_{\mathcal{E}\otimes\mathbb{C}^n} \\
 &=  \sum_{|t|>|t'|}\, \sum_{t(w)=t,\, t(w')=t'}\, \sum_{1\leq i_{|t'|+1},\ldots,i_{m-|t|+|t'|}\leq N+1}\left\langle  K_{w,w'}\right. \\
 &\times  \left. u^{(t')}_{k_1+1,\ldots,k_{|w'|}+1,i_{|t'|+1},\ldots,i_{m-|t|+|t'|},1,\ldots,1}, u^{(t)}_{j_1+1,\ldots,j_{|w|}+1,i_{|t'|+1},\ldots,i_{m-|t|+|t'|}} \right\rangle_{\mathcal{E}} ;
 \end{align*}
\begin{align*}
  C&= \sum_{|t|<|t'|}\left\langle P_{t,t'}u_{t'},u_t\right\rangle_{\mathcal{E}\otimes\mathbb{C}^n}  \\
 &= \sum_{|t|<|t'|}\left\langle \left(\sum_{w,w'\in\free_N:\ t(w)=t,\, t(w')=t'}K_{w,w'}\otimes Z^w{Z^*}^{{w'}^\trans}\right)u_{t'},u_t\right\rangle_{\mathcal{E}\otimes\mathbb{C}^n} \\ 
&= \sum_{|t|<|t'|}\, \sum_{t(w)=t,\, t(w')=t'}\left\langle  K_{w,w'}\otimes \left(E_{j_1+1,1}\otimes\cdots\otimes E_{j_{|w|}+1,1}\otimes I_{N+1}^{\otimes m-|w|}\right)S^{|w|} \right.\\
 &\times   \left. S^{*|w'|}\left(E_{1,k_1+1}\otimes\cdots\otimes E_{1,k_{|w'|}+1}\otimes I_{N+1}^{\otimes m-|w'|}\right) u_{t'},u_t\right\rangle_{\mathcal{E}\otimes\mathbb{C}^n} \\
 &= \sum_{|t|<|t'|}\, \sum_{t(w)=t,\, t(w')=t'}\left\langle  K_{w,w'}\otimes E_{j_1+1,k_{|w'|-|w|+1}+1}\otimes\cdots\otimes E_{j_{|w|}+1,k_{|w'|}+1}\otimes\right. \\ 
 &\otimes  I_{N+1}^{\otimes m-|t'|}\otimes E_{1,k_1+1}\otimes\cdots\otimes E_{1,k_{|w'|-|w|}+1}  \\
 &\times   \sum_{1\leq i_1,\ldots,i_m\leq N+1}u^{(t')}_{i_1,\ldots,i_m}\otimes e_{i_{|w'|-|w|+1}}\otimes\cdots\otimes e_{i_m}\otimes e_{i_1}\otimes\cdots\otimes e_{i_{|w'|-|w|}},\\
 &  \left. \sum_{1\leq s_1,\ldots,s_m\leq N+1}u^{(t)}_{s_1,\ldots,s_m}\otimes e_{s_1}\otimes\cdots\otimes e_{s_m}\right\rangle_{\mathcal{E}\otimes\mathbb{C}^n} \\
 &=  \sum_{|t|<|t'|}\, \sum_{t(w)=t,\, t(w')=t'}\, \sum_{1\leq i_{|t'|+1},\ldots,i_m\leq N+1}\left\langle  K_{w,w'}\right. \\
 &\times  \left. u^{(t')}_{k_1+1,\ldots,k_{|w'|}+1,i_{|t'|+1},\ldots,i_m}, u^{(t)}_{j_1+1,\ldots,j_{|w|}+1,i_{|t'|+1},\ldots,i_m,1\ldots,1} \right\rangle_{\mathcal{E}}.
  \end{align*}

Now, for $t\in\mathbb{Z}^N_+:\ 0\leq |t|\leq m$, 
and $w=g_{j_1}\cdots g_{j_{|w|}}\in\free_N:\ t(w)=t$, 
set 
$$
u^{(t)}_{j_1+1,\ldots,j_{|w|}+1,i_{|t|+1},\ldots,i_m}=0\quad  
\mbox{for}\ (i_{|t|+1},\ldots,i_m)\neq (1,\ldots,1),
$$
and
$$
u^{(t)}_{j_1+1,\ldots,j_{|w|}+1,1,\ldots,1}=
u^{(t(w))}_{j_1+1,\ldots,j_{|w|}+1,1,\ldots,1}=h_w\in\mathcal{E}.
$$
Then the inequality $A+B+C\geq 0$ turns into
$$
\sum_{w,w'\in\free_N:\ 0\leq |w|,|w'|\leq m}
\left\langle K_{w,w'}h_{w'},h_w\right\rangle_{\mathcal{E}}\geq 0.
$$
Since the integer $m>1$ and the vectors $h_w\in\mathcal{E}$ are chosen arbitrarily, 
this means that $K(z,z')$ is a positive non-commutative kernel (for the cases $m=0$ and $m=1$ the inequality above  follows now automatically).
\end{proof}

\section{Proof of Theorem~\ref{conv}}
\label{s:proof2}

For the proof of Theorem~\ref{conv} we need the following lemma. (For a result of a similar nature in a commutative setting, see \cite[Proposition 2.4]{FRKHS}.)
\begin{lemma}\label{l:conv}
An arbitrary convergent positive non-commutative kernel $K(z,z') \in {\mathcal L}({\mathcal E})\langle\langle z, z' \rangle\rangle$ admits a factorization $K(z,z')=H(z)H(z')^*$ with a convergent formal power series $H(z)\in {\mathcal L}(\mathcal{H,E})\langle\langle z, z' \rangle\rangle$, for an auxiliary Hilbert space $\mathcal{H}$. (A formal power series $H(z)\in {\mathcal L}(\mathcal{H,E})\langle\langle z, z' \rangle\rangle$ is called convergent if 
for every $n\in\mathbb{N}$ there exists a connected open neighbourhood $V_n$ of $0$ in
$(\mathbb{C}^{n\times n})^N$ such that the series $H(Z)$ as in \eqref{zpssub}
converges uniformly on compact subsets of $V_n$ in the norm of the Banach space $\mathcal{L}(\mathcal{H} \otimes {\mathbb C}^n,\mathcal{E}\otimes {\mathbb C}^n)$.)
\end{lemma}
\begin{proof}
Let $(\lambda',\lambda'')\in U_1\times U_1$, where $U_n$ is a neighbourhood of $0$ in $(\mathbb{C}^{n\times n})^N$ associated with $K(z,z')$ as in the definition of a convergent non-commutative formal power series which was given in Section~\ref{s:intro}. Then the series $K(\lambda',\lambda'')=\sum_{w\in\mathcal{F}_N}K_{w,w'}{\lambda'}^{t(w)}{\lambda''}^{t(w')}$ converges in the operator norm, and hence there exists a constant $M>0$ such that $\| K_{w,w'}\| {\rho'}^{t(w)}{\rho''}^{t(w')}\leq M$ for all $w,w'\in\mathcal{F}_N$,
where $\rho'=(|\lambda'_1|,\ldots,|\lambda'_N|),\ \rho''=(|\lambda''_1|,\ldots,|\lambda''_N|)$, and for a $w\in\mathcal{F}_N$ the $N$-tuple of integers $t(w)$ is defined in Section~\ref{s:proof1}.   Set $\rho_k=\min\,(|\lambda'_k|,|\lambda''_k|),\ k=1,\ldots,N$, and $\rho=(\rho_1,\ldots,\rho_N)$. Since $K(z,z')$ is a non-commutative positive kernel, it  admits a factorization $K(z,z')=H(z)H(z')^*$ with a  formal power series $H(z)\in {\mathcal L}(\mathcal{H,E})\langle\langle z, z' \rangle\rangle$, for an auxiliary Hilbert space $\mathcal{H}$. This means that $K_{w,w'}=H_wH_{w'}^*$ for every $w,w'\in\mathcal{F}_N$. Therefore we have $$\| K_{w,w}\| {\rho}^{2t(w)}=\| H_wH_w^*\|{\rho}^{2t(w)}=\| H_w\|^2{\rho}^{2t(w)}\leq M,\quad w\in\mathcal{F}_N.$$
 We will use now an argument similar to the one used in the proof of Abel's lemma in the setting of several complex variables (see, e.g., \cite{Sh}). Denote by $V_n$ the open neighbourhood of $0$ in $(\mathbb{C}^{n\times n})^N$ defined by the inequality $\sum_{k=1}^N\rho_k^{-1}\| Z_k\| <1$. Then for $Z\in V_n$ we have:
\begin{eqnarray*}
\lefteqn{\sum_{w\in\mathcal{F}_N}\| H_w\otimes Z^w\| = \sum_{w\in\mathcal{F}_N}\| H_w\|\cdot\| Z^w\|\leq \sum_{w\in\mathcal{F}_N}\| H_w\|\prod_{k=1}^N\| Z_k\|^{t_k(w)} }\\
&\leq & \sqrt{M}\sum_{w\in\mathcal{F}_N}\prod_{k=1}^N(\rho_k^{-1}\| Z_k\| )^{t_k(w)}= \sqrt{M}\sum_{m=0}^\infty\left(\sum_{k=1}^N\rho_k^{-1}\| Z_k\|\right)^m\\
&=& \sqrt{M}\left(1-\sum_{k=1}^N\rho_k^{-1}\| Z_k\|\right)^{-1}<\infty.
\end{eqnarray*}
Therefore for every $n\in\mathbb{N}$ the power series $H(Z)=\sum_{w\in\mathcal{F}_N} H_w\otimes Z^w$ converges absolutely and uniformly on compact subsets of $V_n$ in the norm of ${\mathcal L}(\mathcal{H}\otimes\mathbb{C}^n,\mathcal{E}\otimes\mathbb{C}^n)$, i.e., $H(z)$ is a convergent formal power series.
\end{proof}
\begin{remark}\label{r:nbhd}
It is clear from the proof of Lemma~\ref{l:conv} that the assumption of convergence of the series $K(\lambda',\lambda'')$ at only one point $(\lambda',\lambda'')\in\mathbb{C}^N\times\mathbb{C}^N$ is already sufficient for the conclusion of this lemma. Moreover, this guarantees the existence of neighbourhoods $V_n\ (n\in\mathbb{N})$ ``of the same size", i.e., with the same bounds on the norms of $Z_k,\ k=1,\ldots,N$, such that the series $H(Z)$ converges uniformly and absolutely on compact subsets of $V_n$. 
\end{remark}
\begin{proof}[Proof of Theorem~\ref{conv}.]
 Assume that the convergent formal power series $K(z,z') \in {\mathcal L}({\mathcal E})\langle\langle z, z' \rangle\rangle$ (with respect to a set of neighbourhoods $U_n$ of $0$ in $(\mathbb{C}^{n\times n})^N,\ n\in\mathbb{N}$) is a positive non-commutative kernel. Then by Lemma~\ref{l:conv} $K(z,z')$ admits a factorization
$K(z,z') = H(z) H(z')^*$ for some convergent formal power series
$H(z)$ as in \eqref{zps}, i.e.,
 for every $n\in\mathbb{N}$ there exists a neighbourhood $V_n$ of $0$ in
$(\mathbb{C}^{n\times n})^N$,  such that the series
$H(Z) = \sum_{w \in \free_N} H_w \otimes Z^w$
converges uniformly on compact subsets of $V_n$ (in the norm of the Banach space 
${\mathcal L}({\mathcal H} \otimes {\mathbb C}^n,{\mathcal E} \otimes {\mathbb C}^n)$).
Since we can replace $V_n$ by $V_n\cap U_n$ for all  $n\in\mathbb{N}$, without loss of generality we may assume that $V_n\subset U_n,\ n\in\mathbb{N}$.
In this case $K(Z,Z')$ is a positive
${\mathcal L}({\mathcal E} \otimes {\mathbb C}^n)$-valued kernel on $U_n$.
(We first see that $K(Z,Z') = H(Z) H(Z')^*$ is a positive kernel on $V_n$, then observe that due to uniform convergence of the series for $K(Z,Z')$ on compact subsets of $U_n$ the function $K(Z,Z')$ is sesquianalytic,
and then use the fact that for sesquianalytic kernels the positivity property
propagates to an arbitrarily large region on which the kernel is sesquianalytic;
see \cite[Theorem 1.1.4]{ADRdS} in the one dimensional case -- the general case
is completely analogous.)

Conversely, if
$K(z,z') \in {\mathcal L}({\mathcal E})\langle \langle z, z' \rangle\rangle$
is a convergent formal power series,
and if $K(Z,Z')$
is a positive ${\mathcal L}(\mathcal{E}\otimes\mathbb{C}^n)$-valued
kernel on $U_n$ for every $n\in\mathbb{N}$,
then $K(z,z')$ is a positive non-commutative kernel.
The proof is the same as in Section~\ref{s:proof1}, except that the polynomial
 $P(\lambda,\lambda')$ is guaranteed to be a
positive kernel only on a neighbourhood of $0$ in ${\mathbb C}^N$;
it is clear from the proof of Lemma \ref{lem} that this still implies positive semidefiniteness of the block operator
matrix $M_P$
(alternatively, one can use the propagation property of positivity for sesquianalytic kernels,
mentioned above, to show first that $P(\lambda,\lambda')$ is a positive kernel
on all of ${\mathbb C}^N$).
\end{proof}

\section{Proof of Theorems~\ref{impr}~and~\ref{fact}}
\label{s:proofs3,4}
\begin{proof}[Proof of Theorem \ref{impr}.]
As we already said in Section~\ref{s:intro}, one direction of this theorem is obvious, and thus only the other direction is left to prove. Assume that a formal power series $K(z,z')\in\mathcal{L(E)}\left\langle \left\langle z,z'\right\rangle\right\rangle$ as in \eqref{fps} is such that for every $n\in\mathbb{N}$ and $Z\in {\rm Nilp}_N(n)$ the operator $K(Z,Z)$ as in \eqref{diag} is positive semidefinite. Fix $m\in\mathbb{N}$, and for arbitrary $s>0$ denote by $\mathcal{H}_s$ the Hilbert space with orthogonal basis equal to the set $\mathcal{F}_N^{(m)}$ of words $w\in\mathcal{F}_N$ of length at most $m$, normalized by $\left\langle w,w\right\rangle =s^{-|w|}$. Let $S^*$ denote the backward shift $N$-tuple on $\mathcal{H}_s$ defined as $S_j^*w=v$ if $w=g_jv$, and $S_j^*w=0$ otherwise. Clearly, $S\in {\rm Nilp}_N(\sum_{j=0}^mN^j,m+1)$. Given vectors $h_\alpha\in\mathcal{E}\ (\alpha\in\mathcal{F}_N^{(m)})$, compute:
\allowdisplaybreaks
\begin{align*}
&\left\langle  K(S,S)\sum_{\alpha\in\mathcal{F}_N^{(m)}}h_\alpha\otimes\alpha,\sum_{\beta\in\mathcal{F}_N^{(m)}}h_\beta\otimes\beta
\right\rangle_{\mathcal{E}\otimes\mathcal{H}_s} \\
&= \left\langle \left(\sum_{w,w'\in\mathcal{F}_N^{(m)}}K_{w,w'}\otimes S^wS^{*{w'}^T}\right)\sum_{\alpha\in\mathcal{F}_N^{(m)}}h_\alpha\otimes\alpha,\sum_{\beta\in\mathcal{F}_N^{(m)}}h_\beta\otimes\beta\right\rangle_{\mathcal{E}\otimes\mathcal{H}_s} \\
&= \sum_{w,w'\in\mathcal{F}_N^{(m)}}\sum_{\alpha,\beta\in\mathcal{F}_N^{(m)}}\left\langle K_{w,w'}h_\alpha,h_\beta\right\rangle_\mathcal{E}\left\langle S^{*{w'}^T}\alpha,S^{*w^T}\beta\right\rangle_{\mathcal{H}_s}\\
&= \sum_{\gamma\in\mathcal{F}_N^{(m)}}\sum_{w,w'\in\mathcal{F}_N^{(m-|\gamma |)}}\left\langle K_{w,w'}h_{w'\gamma},h_{w\gamma}\right\rangle_\mathcal{E}\left\langle \gamma,\gamma\right\rangle_{\mathcal{H}_s}\\
&= \sum_{\gamma\in\mathcal{F}_N^{(m)}}\sum_{w,w'\in\mathcal{F}_N^{(m-|\gamma |)}}\left\langle K_{w,w'}h_{w'\gamma},h_{w\gamma}\right\rangle_\mathcal{E}s^{-|\gamma |}.
\end{align*}
By hypothesis, $K(S,S)$ is a positive semidefinite operator on $\mathcal{H}_s$ for each $s>0$. Letting $s$ tend to infinity gives
$$  \sum_{w,w'\in\mathcal{F}_N^{(m)}}\left\langle K_{w,w'}h_{w'},h_w\right\rangle_\mathcal{E}\geq 0,$$
which, in virtue of arbitrariness of $m\in\mathbb{N}$ and of vectors $h_w\in\mathcal{E}\ (w\in\mathcal{F}_N^{(m)})$, means that $K(z,z')$ is a positive non-commutative kernel.
\end{proof}
\begin{proof}[Proof of Theorem~\ref{fact}.]
The statement follows from Theorem~\ref{impr}, since in this case $K(z,z')\in\mathcal{L(E)}\left\langle z,z'\right\rangle$ is a positive non-commutative kernel. The estimate for $\dim (\mathcal{H})$ is easily obtained from the fact that the  matrix $(K_{w,w'})_{|w|\leq m,|w'|\leq m}$ with blocks from $\mathcal{L(E)}$ is of size $\left(\sum_{j=0}^mN^j\right)\times\left(\sum_{j=0}^mN^j\right)$ and positive semidefinite.
The last statement of the theorem (about the size and the joint rank of nilpotency 
of the test matrices) follows from the proof of Theorem~\ref{impr}.
\end{proof}

\providecommand{\bysame}{\leavevmode\hbox to3em{\hrulefill}\thinspace}
\providecommand{\MR}{\relax\ifhmode\unskip\space\fi MR }
\providecommand{\MRhref}[2]{%
  \href{http://www.ams.org/mathscinet-getitem?mr=#1}{#2}
}
\providecommand{\href}[2]{#2}

\end{document}